\documentstyle[leqno]{article}

\newtheorem{th}{Theorem}[section]
\newtheorem{prop}[th]{Proposition}
\newtheorem{cor}[th]{Corollary}
\newcounter{defin}[section]
\renewcommand{\thedefin}{\thesection.\arabic{defin}}
\newcounter{ex}[section]
\renewcommand{\theex}{\thesection.\arabic{ex}}
\newcounter{rem}[section]
\renewcommand{\therem}{\thesection.\arabic{rem}}
\title{The Geometry of Metrical Multi-Time \\
Lagrange Spaces}
\author{Mircea Neagu and Constantin Udri\c ste}
\date{}
\begin{document}
\maketitle
\begin{abstract}
Section 1 contains historical and bibliographical notes upon the problem of
geometrization of Lagrangians defined on the tangent bundle or the jet bundle
of order one, and  emphasizes the original elements of our approach in this
direction. The geometrization of a  Kronecker $h$-regular Lagrangian function
with partial
derivatives begins in Section 2 by introduction of notion of metrical
multi-time Lagrange space $ML^n_p=(J^1(T,M),L)$  and by proving a
theorem of characterization of these spaces. Section 3 constructs the
canonical nonlinear connection $\Gamma=(M^{(i)}_{(\alpha)\beta}, N^{(i)}_
{(\alpha)j})$, naturally induced by the Lagrangian ${\cal L}=L\sqrt{\vert h\vert}$
of the metrical multi-time Lagrange space $ML^n_p$. At the same time, Section
3 offers a geometrical interpretation  to the extremals of the
Lagrangian ${\cal L}$. Section 4 proves the theorem of existence  and
uniqueness of Cartan canonical connection $C\Gamma$ of a metrical multi-time
Lagrange space $ML^n_p$ and studies its torsion and curvature d-tensors.
\end{abstract}
{\bf Mathematics Subject Classification (2000):} 53B05, 53B40, 53C43.\\
{\bf Key  words:} temporal and spatial sprays, harmonic maps, metrical
multi-time Lagrange space, canonical Cartan metrical connection, torsion and
curvature d-tensors.

\section{Introduction}

\hspace{5mm} A lot of geometrical models in Mechanics, Physics or Biology are
based on the notion of ordinary Lagrangian. In this sense, we recall that  a Lagrange  space
$L^n=(M,L(x,y))$ is defined as a pair which consists of a real, smooth,\linebreak
$n$-dimensional manifold $M$ coordinated by $(x^i)_{i=\overline{1,n}}$ and a regular Lagrangian
$L:TM\to R$. The differential geometry of Lagrange spaces is now
used in various fields to study natural phenomena
where the dependence on position, velocity  or momentum is involved \cite{7}.
Also, this geometry gives a model for both the gravitational
and electromagnetic field theory, in a very natural blending of the geometrical
structure of the space with the characteristic properties of the physical
fields.

At the same time, there are many problems in Physics and  Variational Calculus
in which time dependent Lagrangians (i. e., a  smooth real function on
$R\times TM$) are involved. A geometrization of time dependent Lagrangians
was realized by Miron and Anastasiei in \cite{7}, using the
configuration bundle $R\times TM\to M$, whose geometrical invariance group
(gauge group) is of the form
\begin{equation}\label{G_1}
\left\{\begin{array}{l}\medskip
\tilde t=t\\\medskip
\tilde x^i=\tilde x^i(x^j)\\
\displaystyle{\tilde y^i={\partial\tilde x^i\over\partial x^j}y^j}.
\end{array}\right.
\end{equation}
The attached geometry is called the {\it "Rheonomic Lagrange Geometry"}.
Nevertheless, the main inconvenient of rheonomic Lagrange geometry is
determined by the {\it "absolute"} character of time $t$ which is emphasized
by the structure of its gauge group.
In our paper, we remove this inconvenient, replacing
the above bundle of configuration with 1-jet fibre bundle $J^1(R,M)\equiv R\times TM\to
R\times M$ which is characterized by the gauge group
\begin{equation}\label{G_2}
\left\{\begin{array}{l}\medskip
\tilde t=\tilde t(t)\\\medskip
\tilde x^i=\tilde x^i(x^j)\\
\displaystyle{\tilde y^i={\partial\tilde x^i\over\partial x^j}
{dt\over d\tilde t}y^j}.
\end{array}\right.
\end{equation}
The structure of our gauge group underlines the {\it relativistic}
role of time,  and consequently, we can use  the name of {\it "Relativistic Rheonomic
Lagrange Geometry"}.

In the last thirty years, many mathematiciens were concerned by  the
geometrization of a {\it multi-time Lagrangian} depending on first order
partial derivatives, which is defined on the 1-jet fibre bundle $J^1(T,M)$,
where $T$ is a smooth, real, $p$-dimensional {\it "multi-time"} manifold
coordinated by $(t^\alpha)_{\alpha=\overline{1,p}}$ and $M$ is a smooth,  real,
$n$-dimensional {\it "spatial"} manifold coordinated by $(x^i)_{i=\overline{1,n}}$.

One point of view is described by Gotay, Isenberg and Marsden \cite{4} and
is known under the name of {\it Multisimplectic Geometry}. This paper
naturally generalizes the dual Hamiltonian formalism used  in the
classical mechanics and stands out by its {\it finite dimensional} and {\it
non metric} spatial model $J^1(T,M)$.
In contrast, Michor and Ra\c tiu \cite{6} construct their geometrization on the
{\it infinite dimensional} space of the embeddings $Emb(T,M)$ which is
endowed with a {\it metrical structure} $G$.
The third way is sketched by Miron, Kirkovits and Anastasiei
in \cite{8}. They construct a {\it metrical} geometry attached
to a first order Lagrangian on the {\it finite dimensional} total space of the
vector bundle $\oplus^p_{\alpha=1}\;TM\to M$, where  the coordinates of
$\alpha$-th copy of $TM$ are denoted $(x^i,x^i_\alpha)$, and the gauge group of
bundle of configuration is of the form
\begin{equation}
\left\{\begin{array}{l}\medskip
\tilde x^i=\tilde x^i(x^j)\\
\displaystyle{\tilde x^i_\alpha={\partial\tilde x^i\over\partial x^j}x^j_\alpha}.
\end{array}\right.
\end{equation}
This geometry relies on a given semi-Riemannian
metric $h_{\alpha\beta}$ on $R^p$ and a given {\it "a priori"} nonlinear connection
on $E=\oplus^p_{\alpha=1}\;TM$. Starting with a multi-time Lagrangian
${\cal L}:\oplus^p_{\alpha=1}\;TM\to R$  and using an adapted basis of the nonlinear
connection, they introduce a Sasakian-like metric on $TE$, setting
\begin{equation}
G=g_{ij}\delta x^i\otimes\delta x^j+G^{(\alpha)(\beta)}_{(i)(j)}\delta x^i_\alpha
\otimes\delta x^j_\beta,
\end{equation}
where $\displaystyle{G^{(\alpha)(\beta)}_{(i)(j)}(x^k,x^k_\gamma)={1\over 2}
{\partial^2{\cal L}\over\partial x^i_\alpha\partial x^j_\beta}}$ and $g_{ij}
(x^k,x^k_\gamma)=h_{\alpha\beta}G^{(\alpha)(\beta)}_{(i)(j)}$. Also,
the paper \cite{8} use  the {\it Lagrangian
density}
$
{\cal D}={\cal L}\;dt^1\wedge dt^2\ldots\wedge dt^p,
$
developing a multi-time Lagrangian geometry, in
the sense of linear connections, torsions and curvatures.

In this paper, we naturally extend the Rheonomic Lagrange Geometry
of {\it vector} bundle $R\times TM\to M$ to the 1-jet fibre bundle
$J^1(T,M)\to T\times M$.  Using  the gauge group
\begin{equation}
\left\{\begin{array}{l}\medskip
\tilde t^\alpha=\tilde t^\alpha(t^\beta)\\\medskip
\tilde x^i=\tilde x^i(x^j)\\
\displaystyle{\tilde x^i_\alpha={\partial\tilde x^i\over\partial x^j}
{\partial t^\beta\over\partial\tilde t^\alpha}x^j_\beta},
\end{array}\right.
\end{equation}
which is more general than that used in the papers  \cite{7}, \cite{8}. In order
to develope our geometry, we start {\it "a priori"} with a semi-Riemannian
metric $h=h_{\alpha\beta}(t)$ on the temporal manifold $T$, and we use the
following three distinct notions:

i) {\it multi-time Lagrangian function} $-$ A smooth function $L:J^1(T,M)\to R$.

ii) {\it multi-time Lagrangian} (Olver's terminology) $-$ A local function ${\cal L}$
on $J^1(T,M)$ which transform by the rule $\tilde{\cal L}={\cal L}\vert\det J\vert$, where
$J$ is the Jacobian matrix of coordinate transformations $t^\alpha=
t^\alpha(\tilde t^\beta)$. If $L$ is a Lagrangian function on
1-jet fibre bundle, then ${\cal L}=L\sqrt{\vert h\vert}$
represent a Lagrangian on $J^1(T,M)$.

iii) {\it multi-time Lagrangian density} (Marsden's terminology) $-$ A smooth map
${\cal D}:J^1(T,M)\to \Lambda^p(T^*T)$. For example, the entity ${\cal D}=
{\cal L}dt^1\wedge dt^2\wedge\ldots\wedge dt^p$, where ${\cal L}$ is a
Lagrangian, represents a Lagrangian density on $J^1(T,M)$.

In this terminology, we create a geometry attached to a first
order {\it Kronecker $h$-regular Lagrangian function} on $J^1(T,M)$, which
can be called the {\it Metrical Multi-Time Lagrangian Geometry}. The condition of
Kronecker $h$-regularity imposed to the given multi-time Lagrangian function
$L$ is
$$
G^{(\alpha)(\beta)}_{(i)(j)}(t^\gamma,x^k,x^k_\gamma)={1\over 2}
{\partial^2L\over\partial x^i_\alpha\partial x^j_\beta}=
h^{\alpha\beta}(t^\gamma)g_{ij}(t^\gamma,x^k,x^k_\gamma),
$$
where $g_{ij}(t^\gamma,x^k,x^k_\gamma)$ is a d-tensor on $J^1(T,M)$, symmetric,
having  the rank $n$, and of constant signature. This condition allows us
to build from $L$ a natural nonlinear connection on $J^1(T,M)$.
At the same time, the condition of Kronecker $h$-regularity  is required by the
following reasons:

(1) the construction of Sasakian-like metric
\begin{equation}
G=h_{\alpha\beta}dt^\alpha\otimes dt^\beta+g_{ij}dx^i\otimes dx^j+
h^{\alpha\beta}g_{ij}\delta x^i_\alpha\otimes\delta x^j_\beta,
\end{equation}
on jet bundle $J^1(T,M)$, that has  the physical meaning of gravitational
potential on $J^1(T,M)$;

(2) to  find explicitly canonical d-connections, d-torsions and d-curvatures.

We emphasize that the physical aspects of the metrical multi-time Lagrange
geometry are exposed in \cite{9},  \cite{10}. From this point of view, the Maxwell and
Einstein equations allow us to appreciate
the Metrical Multi-Time Lagrange Geometry like a natural model
necessary in the study of physical fields in a general setting.

Finally, we point out again that the Lagrangian density used in our study is
\begin{equation}
{\cal D}=L\sqrt{\vert h\vert}\;dt^1\wedge dt^2\ldots\wedge dt^p.
\end{equation}

\section{Metrical multi-time Lagrange spaces}

\setcounter{equation}{0}
\hspace{5mm}  Let us consider $T$ (resp. $M$) a {\it "temporal"}  (resp. {\it
"spatial"}) manifold of dimension $p$ (resp. $n$), coordinated by $(t^\alpha)_
{\alpha=\overline{1,p}}$ (resp. $(x^i)_{i=\overline{1,n}}$). Let $E=J^1(T,M)\to
T\times M$ be the jet fibre bundle of order one associated to these manifolds.
The {\it bundle of configuration} $J^1(T,M)$ is coordinated by
$(t^\alpha,x^i,x^i_\alpha)$, where $\alpha=\overline{1,p}$ and $i=\overline{1,n}$.
Note that the terminology used above is justified in \cite{12}.\medskip\\
\addtocounter{rem}{1}
{\bf Remarks \therem} i) Throughout  this paper, the indices $\alpha,\beta,\gamma,
\ldots$ run from $1$ to $p$, and the indices $i,j,k,\ldots$ run from $1$ to $n$.

ii) In the particular case $T=R$ (i. e., the temporal manifold $T$ is the usual
time axis represented by the set of real numbers), the coordinates  $(t^1,x^i,
x^i_1)$ of the 1-jet space  $J^1(R,M)\equiv R\times TM$ are denoted $(t,x^i,y^i)$.
\medskip

We start our study considering  a smooth multi-time Lagrangian function
$L:E\to R$, which is locally expressed by $E\ni(t^\alpha,x^i,x^i_\alpha)\to
L(t^\alpha,x^i,x^i_\alpha)\in R$. The {\it vertical fundamental  metrical
d-tensor} of $L$ is
\begin{equation}
G^{(\alpha)(\beta)}_{(i)(j)}={1\over 2}{\partial^2L\over\partial x^i_\alpha
\partial x^j_\beta}.
\end{equation}

Now, let $h=(h_{\alpha\beta})$ be a fixed semi-Riemannian metric on the temporal
manifold $T$ and $g_{ij}(t^\gamma, x^k, x^k_\gamma)$ be a d-tensor on $E$, symmetric,
of rank $n$, and having a constant signature.\medskip\\
\addtocounter{defin}{1}
{\bf Definition \thedefin} A multi-time Lagrangian function $L:E\to R$ whose vertical
fundamental metrical d-tensor is of the form
\begin{equation}
G^{(\alpha)(\beta)}_{(i)(j)}(t^\gamma,x^k,x^k_\gamma)=h^{\alpha\beta}(t^
\gamma)g_{ij}(t^\gamma,x^k,x^k_\gamma),
\end{equation}
is called a {\it Kronecker $h$-regular multi-time Lagrangian function
with respect to the
temporal semi-Riemannian metric $h=(h_{\alpha\beta})$}.\medskip

In this context, we can introduce the  following\medskip\\
\addtocounter{defin}{1}
{\bf  Definition \thedefin} A pair $ML^n_p=(J^1(T,M),L)$, where $p=\dim T$ and $n=\dim M$,
which consists of the  1-jet  fibre bundle and a
Kronecker $h$-regular multi-time Lagrangian function $L:J^1(T,M)\to R$ is called a {\it
metrical multi-time Lagrange space}.\medskip\\
\addtocounter{rem}{1}
{\bf Remarks \therem} i) In the particular case $(T,h)=(R,\delta)$, a metrical
multi-time Lagrange space is called a {\it  relativistic rheonomic Lagrange
space}  and is denoted $RL^n=(J^1(R,M),L)$.

ii) If  the temporal manifold $T$ is 1-dimensional, then,
via a temporal reparametrization, we have  $J^1(T,M)\equiv J^1(R,M)$.
In other words, a metrical multi-time Lagrangian space having $\dim T=1$ is a
{\it reparametrized relativistic rheonomic Lagrange space}.\medskip\\
\addtocounter{ex}{1}
{\bf Examples \theex} i) Suppose that the spatial manifold $M$ is also endowed with
a semi-Riemannian metric $g=(g_{ij}(x))$. Then, the  multi-time Lagrangian function
\begin{equation}
L_1:J^1(T,M)\to R,\quad
L_1=h^{\alpha\beta}(t)g_{ij}(x)x^i_\alpha x^j_\beta
\end{equation}
is  a Kronecker $h$-regular multi-time Lagrangian function. Consequently,
$ML^n_p=(J^1(T,M),L_1)$ is
a metrical multi-time Lagrange space. We underline that the multi-time Lagrangian
${\cal L}_1=L_1\sqrt{\vert h\vert}$ is exactly the energy multi-time Lagrangian
whose extremals are the harmonic maps between the pseudo-Riemannian manifolds
$(T,h)$ and $(M,g)$ \cite{3}. At the same time, this
multi-time Lagrangian is a basic object in the physical theory of bosonic
strings.

ii) In above notations, taking $U^{(\alpha)}_{(i)}(t,x)$ as a d-tensor field on
$E$ and \linebreak $F:T\times M\to R$ a smooth map, the more general multi-time
Lagrangian function
\begin{equation}
L_2:E\to R,\quad
L_2=h^{\alpha\beta}(t)g_{ij}(x)x^i_\alpha x^j_\beta+U^{(\alpha)}_{(i)}(t,x)x^i_
\alpha+F(t,x)
\end{equation}
is also a Kronecker $h$-regular multi-time Lagrangian. The metrical multi-time
Lagrange space $ML^n_p=(J^1(T,M),L_2)$ is called the {\it autonomous metrical
multi-time Lagrange space  of electrodynamics} because, in the particular case
$(T,h)=(R,\delta)$, we recover the classical Lagrangian space of electrodynamics
\cite{7} which governs the movement law  of a particle placed concomitantly into
a gravitational field and an electromagnetic one. From physical point of view,
the semi-Riemannian metric $h_{\alpha\beta}(t)$ (resp. $g_{ij}(x)$) represents the {\it gravitational
potentials} of the space $T$ (resp. $M$), the d-tensor $U^{(\alpha)}_{(i)}(t,x)$
stands for the {\it electromagnetic potentials} and $F$ is a function which is
called {\it potential function}. The  non-dynamical character of spatial gravitational
potentials $g_{ij}(x)$ motivates us to use the term {\it "autonomous"}.

iii) More general, if  we consider $g_{ij}(t,x)$ a d-tensor field  on $E$,
symmetric,  of rank $n$ and  having constant signature on $E$, we can
define the Kronecker $h$-regular multi-time Lagrangian function
\begin{equation}
L_3:E\to R,\quad
L_3=h^{\alpha\beta}(t)g_{ij}(t,x)x^i_\alpha x^j_\beta+U^{(\alpha)}_{(i)}(t,x)x^i_
\alpha+F(t,x).
\end{equation}
The pair $ML^n_p=(J^1(T,M),L_3)$ is a metrical multi-time Lagrange  space which
is called the
{\it non-autonomous metrical multi-time Lagrange space of electrodynamics}.
Physically, we
remark that the gravitational potentials $g_{ij}(t,x)$ of the spatial manifold
$M$ are dependent of the temporal coordinates $t^\gamma$, emphasizing their
dynamic character.

An important role and, at the same time, an obstruction in the subsequent
development  of the metrical multi-time Lagrangian theory, is played by  the following
\begin{th} {(characterization of metrical multi-time Lagrange spaces)}\\
If we have $\dim T\geq 2$, then the following statements are equivalent:

i) $L$ is a Kronecker $h$-regular multi-time Lagrangian function on $J^1(T,M)$.

ii) The multi-time Lagrangian function $L$ reduces  to a non-autonomous electrodynamics
Lagrangian function, that is,
$$
L=h^{\alpha\beta}(t)g_{ij}(t,x)x^i_\alpha x^j_\beta+U^{(\alpha)}_{(i)}(t,x)x^i_
\alpha+F(t,x).
$$
\end{th}
{\bf Proof.} ii)$\Rightarrow$ i) Obviously.\\
\hspace*{12mm} i)$\Rightarrow$ ii) Suppose that $L$ is a Kronecker $h$-regular
multi-time Lagrangian function,  that is,
$$
{1\over 2}{\partial^2L\over\partial x^i_\alpha\partial x^j_\beta}=h^{\alpha\beta}(t^
\gamma)g_{ij}(t^\gamma,x^k,x^k_\gamma).
$$

Firstly, we assume that there are two distinct indices $\alpha$ and $\beta$
in the set $\{1,\ldots,p\}$ such that $h^{\alpha\beta}\neq 0$. Let $k$ (resp.
$\gamma$) be an  arbitrary element of the set $\{1,\ldots,n\}$ (resp. $\{1,
\ldots,p\}$). Differentiating the above relation by $x^k_\gamma$ and using the Schwartz
theorem, we obtain the equalities
$$
{\partial g_{ij}\over\partial x^k_\gamma}h^{\alpha\beta}=
{\partial g_{jk}\over\partial x^i_\alpha}h^{\beta\gamma}=
{\partial g_{ik}\over\partial x^j_\beta}h^{\gamma\alpha},\quad
\forall\; \alpha,\beta,\gamma\in\{1,\ldots,p\},\quad\forall\; i,j,k\in\{1,\ldots,n\}.
$$
Contracting by $h_{\gamma\mu}$, we deduce
$$
{\partial g_{ij}\over\partial x^k_\gamma}h^{\alpha\beta}h_{\gamma\mu}=0,\quad
\forall\; \mu\in\{1,\ldots,p\}.
$$
The assumption $h^{\alpha\beta}\neq 0$ implies that $\displaystyle{{\partial
g_{ij}\over\partial x^k_\gamma}=0}$ for arbitrary $k$ and $\gamma$, that is,
\linebreak $g_{ij}=g_{ij}(t^\mu,x^m)$.

Secondly, assuming that $h^{\alpha\beta}=0,\;\forall\;\alpha\neq\beta\in\{1,
\ldots,p\}$, it follows \linebreak
$h^{\alpha\beta}=h^\alpha\delta^\alpha_\beta,\;
\forall\;\alpha,\beta\in\{1,\ldots,p\}$. In other words, on $T$ we use an
orthogonal system of coordinates. In these conditions, the relations
$$
\begin{array}{l}\medskip
\displaystyle{
{\partial^2L\over\partial x^i_\alpha\partial x^j_\beta}=0,\quad
\forall\;\alpha\neq\beta\in\{1,\ldots,p\},\quad\forall\;i,j\in\{1,\ldots,n\}},
\\
\displaystyle{{1\over 2h^\alpha(t)}{\partial^2L\over\partial x^i_\alpha\partial
x^j_\alpha}=g_{ij}(t^\mu,x^m,x^m_\mu),\quad\forall\;\alpha\in\{1,\ldots,p\}
,\quad\forall\;i,j\in\{1,\ldots,n\}}
\end{array}
$$
are true. Now, if  we fixe the indice  $\alpha$ in the set $\{1,\ldots,p\}$,
we deduce by first relation that the local functions $\displaystyle{{\partial
L\over\partial x^i_\alpha}}$ depend just of the coordinates $(t^\mu,x^m,x^m_
\alpha)$. Considering $\beta\neq\alpha$ in the set $\{1,\ldots,
p\}$, the second relation implies
$$
\displaystyle{{1\over 2h^\alpha(t)}{\partial^2L\over\partial x^i_\alpha\partial
x^j_\alpha}={1\over 2h^\beta(t)}{\partial^2L\over\partial x^i_\beta\partial
x^j_\beta}=g_{ij}(t^\mu,x^m,x^m_\mu),\quad\forall\;i,j\in\{1,\ldots,n\}}.
$$
Because the first term of the above equality depends just of the coordinates
$(t^\mu,x^m,x^m_\alpha)$ while the second term is dependent just of $(t^\mu,
x^m,x^m_\beta)$ and \linebreak
$\alpha\neq\beta$, we conclude that
$g_{ij}=g_{ij}(t^\mu,x^m)$.

Finally, the relation
$$
{1\over 2}{\partial^2L\over\partial x^i_\alpha\partial x^j_\beta}=h^{\alpha\beta}(t^
\gamma)g_{ij}(t^\gamma,x^k),\quad\forall\;\alpha,\beta\in\{1,\ldots,p\},\quad
\forall\;i,j\in\{1,\ldots,n\}
$$
implies without difficulties that the multi-time Lagrangian function  $L$ is
a  non-autonomous multi-time Lagrangian function of electrodynamics. \rule{5pt}{5pt}
\begin{cor}
The vertical fundamental metrical d-tensor of an arbitrary Kronecker
$h$-regular multi-time Lagrangian function $L$ is of the form
\begin{equation}
G^{(\alpha)(\beta)}_{(i)(j)}={1\over 2}{\partial^2L\over\partial x^i_\alpha
\partial x^j_\beta}=\left\{\begin{array}{ll}\medskip
h^{11}(t)g_{ij}(t,x^k,y^k),&\dim T=1\\
h^{\alpha\beta}(t^\gamma)g_{ij}(t^\gamma,x^k),&\dim T\geq 2.
\end{array}\right.
\end{equation}
\end{cor}
\addtocounter{rem}{1}
{\bf Remarks \therem} i) It is obvious that the preceding theorem is  an obstruction
in the development  of a fertile metrical multi-time Lagrangian geometry. This
obstruction will be removed in a subsequent paper by the introduction of a more general notion,
that of {\it generalized metrical multi-time Lagrange space} \cite{9}. The generalized
metrical multi-time Lagrange geometry is constructed using just a given
{vertical Kronecker $h$-regular metrical $d$-tensor $G^{(\alpha)(\beta)}_{(i)(j)}$
on the 1-jet space $J^1(T,M)$.

ii) In the case $\dim T\geq 2$, the above theorem obliges us to continue the
study of the metrical multi-time Lagrangian spaces theory, channeling our
attention upon the non-autonomous electrodynamics metrical multi-time
Lagrangian spaces.

\section{Sprays. Nonlinear connection. Harmonic maps}

\setcounter{equation}{0}
\hspace{5mm} Let $ML^n_p=(J^1(T,M),L)$, where $\dim T=p,\;\dim M=n$, be a metrical
multi-time Lagrange space whose vertical fundamental  metrical d-tensor is
$$
G^{(\alpha)(\beta)}_{(i)(j)}={1\over 2}{\partial^2L\over\partial x^i_\alpha
\partial x^j_\beta}=\left\{\begin{array}{ll}\medskip
h^{11}(t)g_{ij}(t,x^k,y^k),&p=1\\
h^{\alpha\beta}(t^\gamma)g_{ij}(t^\gamma,x^k),&p\geq 2.
\end{array}\right.
$$
Note that all subsequent entities with geometrical or physical meaning will be
directly obtained from the fundamental vertical metrical d-tensor $G^{(\alpha)(\beta)}_
{(i)(j)}$.  This fact points out the {\it metrical character}
and the naturalness of the metrical multi-time Lagrangian geometry that we construct.
At the same  time, the form of the invariance gauge group
$$
\left\{\begin{array}{l}\medskip
\tilde t^\alpha=\tilde t^\alpha(t^\beta)\\\medskip
\tilde x^i=\tilde x^i(x^j)\\
\displaystyle{\tilde x^i_\alpha={\partial\tilde x^i\over\partial x^j}
{\partial t^\beta\over\partial\tilde t^\alpha}x^j_\beta}
\end{array}\right.
$$
of the fibre bundle $J^1(T,M)\to T\times M$ allows us to look out the metrical
multi-time Lagrange geometry as a {\it "parametrized"} theory,  in Marsden's sense
\cite{4}.

Now, assume that the  semi-Riemannian temporal manifold $(T,h)$ is compact
and orientable. In this context, we can define the {\it energy functional} of
the Lagrangian function $L$, setting
$$
{\cal E}_L:C^\infty(T,M)\to R,\quad
{\cal E}_L(f)=\int_TL(t^\alpha,x^i,x^i_\alpha)\sqrt{\vert h\vert}\;dt^1\wedge dt^2
\wedge\ldots\wedge dt^p,
$$
where the smooth map $f$ is locally expressed by $(t^\alpha)\to(x^i(t^\alpha))$
and $\displaystyle{x^i_\alpha={\partial x^i\over\partial t^\alpha}}$.

The extremals of the energy functional ${\cal E}_L$ verifies the Euler-Lagrange
equations for every $i\in\{1,2,\ldots,n\}$,
\begin{equation}
2G^{(\alpha)(\beta)}_{(i)(j)}x^j_{\alpha\beta}+{\partial^2L\over\partial x^j
\partial x^i_\alpha}x^j_\alpha-{\partial L\over\partial x^i}+{\partial^2L\over
\partial t^\alpha\partial x^i_\alpha}+{\partial L\over\partial x^i_\alpha}
H^\gamma_{\alpha\gamma}=0,
\end{equation}
where
$\displaystyle{x^j_{\alpha\beta}={\partial^2x^j\over\partial t^\alpha
\partial t^\beta}}$, and $H^\gamma_{\alpha\beta}$ are the Christoffel
symbols of the semi-Riemannian  metric $h_{\alpha\beta}$.

Taking into account the Kronecker  $h$-regularity of the Lagrangian function $L$, it
is  possible  to rearrange the Euler-Lagrange equations  of Lagrangian
${\cal L}=L\sqrt{\vert h\vert}$ in the form
\begin{equation}
\Delta_hx^k+2{\cal G}^k(t^\mu,x^m,x^m_\mu)=0,
\end{equation}
where
$$
\Delta_hx^k=h^{\alpha\beta}\{x^k_{\alpha\beta}-H^\gamma_{\alpha\beta}
x^k_\gamma\},
$$
$$
2{\cal G}^k={g^{ki}\over 2}\left\{{\partial^2L\over\partial x^j
\partial x^i_\alpha}x^j_\alpha-{\partial L\over\partial x^i}+{\partial^2L\over
\partial t^\alpha\partial x^i_\alpha}+{\partial L\over\partial x^i_\alpha}
H^\gamma_{\alpha\gamma}+2g_{ij}h^{\alpha\beta}H^\gamma_{\alpha\beta}x^j_\gamma
\right\}.
$$

By a  direct calculation, we deduce that the local geometrical entities
\begin{equation}
\begin{array}{l}\medskip
\displaystyle{2{\cal S}^k={g^{ki}\over 2}\left\{{\partial^2L\over\partial x^j
\partial x^i_\alpha}x^j_\alpha-{\partial L\over\partial x^i}\right\}}\\\medskip
\displaystyle{2{\cal H}^k={g^{ki}\over 2}\left\{{\partial^2L\over\partial t^\alpha
\partial x^i_\alpha}+{\partial L\over\partial x^i_\alpha}H^\gamma_{\alpha\gamma}\right\}}\\
2{\cal J}^k=h^{\alpha\beta}H_{\alpha\beta}^\gamma x^j_\gamma
\end{array}
\end{equation}
verify the following transformation rules
\begin{equation}
\begin{array}{l}\medskip
\displaystyle{2{\cal S}^p=2\tilde{\cal S}^r{\partial x^p\over\partial\tilde x^r}+h^{\alpha\mu}
{\partial x^p\over\partial\tilde x^l}{\partial\tilde t^\gamma\over\partial t^\mu}
{\partial\tilde x^l_\gamma\over\partial x^j}x^j_\alpha}\\\medskip
\displaystyle{2{\cal H}^p=2\tilde{\cal H}^r{\partial x^p\over\partial\tilde x^r}+h^{\alpha\mu}
{\partial x^p\over\partial\tilde x^l}{\partial\tilde t^\gamma\over\partial t^\mu}
{\partial\tilde x^l_\gamma\over\partial t^\alpha}}\\
\displaystyle{2{\cal J}^p=2\tilde{\cal J}^r{\partial x^p\over\partial\tilde x^r}-h^{\alpha\mu}
{\partial x^p\over\partial\tilde x^l}{\partial\tilde t^\gamma\over\partial t^\mu}
{\partial\tilde x^l_\gamma\over\partial t^\alpha}}.
\end{array}
\end{equation}
Consequently, the local entities $2{\cal G}^p=2{\cal S}^p+2{\cal H}^p+2{\cal J}^p$
modify by the transformation laws
\begin{equation}
2\tilde{\cal G}^r=2{\cal G}^p{\partial\tilde x^r\over\partial x^p}-h^{\alpha\mu}
{\partial x^p\over\partial\tilde x^j}{\partial\tilde x^r_\mu\over\partial x^p}
\tilde x^j_\alpha.
\end{equation}
and therefore the  geometrical object ${\cal  G}=({\cal G}^r)$ is a
{\it spatial $h$-spray} \cite{12}.

Following the paper \cite{12}, we can offer a  geometrical interpretation to the  equations
$$
\Delta_hx^l+2{\cal G}^l(t^\gamma,x^k,x^k_\gamma)=0,\quad\forall\;l\in\{1,\ldots,
n\},
$$
via the {\it harmonic map equations of a  spatial spray}, if the
spatial $h$-spray ${\cal G}$ is the $h$-trace of a spatial spray $G$.

In the particular case
$\dim T=1$, every spatial $h$-spray ${\cal G}=({\cal G}^l)$ is the $h$-trace
of a spatial spray, namely $G=(G^{(l)}_{(1)1})$, where $G^{(l)}_{(1)1}=h_{11}
{\cal G}^l$. In other words, the equality ${\cal G}^l=h^{11}G^{(l)}_{(1)1}$
is true.

On the other hand, in the case $\dim T\geq 2$, the characterization theorem of
the Kronecker $h$-regular Lagrangians functions ensures us that
$$
L=h^{\alpha\beta}(t)g_{ij}(t,x)x^i_\alpha x^j_\beta+U^{(\alpha)}_{(i)}(t,x)x^i_
\alpha+F(t,x).
$$
In this particular situation, by computations, the  expressions of ${\cal S}^l,\;
{\cal H}^l$ and ${\cal J}^l$ reduce  to
\begin{equation}
\begin{array}{l}\medskip
\displaystyle{2{\cal S}^l=h^{\alpha\beta}\Gamma^l_{jk}x^j_\alpha x^k_\beta+{g^{li}\over 2}
\left[U^{(\alpha)}_{(i)j}x^j_\alpha-{\partial F\over\partial x^i}\right]}
\\\medskip
\displaystyle{2{\cal H}^l=-h^{\alpha\beta}H^\gamma_{\alpha\beta}x^l_\gamma+{g^{li}\over 2}
\left[2h^{\alpha\beta}{\partial g_{ij}\over\partial t^\alpha}x^j_\beta+{\partial
U^{(\alpha)}_{(i)}\over\partial t^\alpha}+U^{(\alpha)}_{(i)}H^\gamma_{\alpha\gamma}
\right]}\\
2{\cal J}^l=h^{\alpha\beta}H^\gamma_{\alpha\beta}x^l_\gamma,
\end{array}
\end{equation}
where
$$
\displaystyle{\Gamma^l_{jk}={g^{li}\over 2}\left({\partial g_{ij}\over\partial
x^k}+{\partial g_{ik}\over\partial x^j}-{\partial g_{jk}\over\partial x^i}
\right)}
$$
are the {\it generalized Christoffel symbols of the "multi-time" dependent
metric $g_{ij}$}, and
$$
\displaystyle{U^{(\alpha)}_{(i)j}={\partial U^{(\alpha)}_{(i)}\over
\partial x^j}-{\partial U^{(\alpha)}_{(j)}\over\partial x^i}}.
$$
Consequently,
the expression of the spatial $h$-spray ${\cal G}=({\cal G}^l)$  becomes
\begin{equation}
2{\cal G}^p=2{\cal S}^p+2{\cal H}^p+2{\cal J}^p=h^{\alpha\beta}\Gamma^l_{jk}
x^j_\alpha x^k_\beta+2{\cal T}^l,
\end{equation}
where
\begin{equation}
2{\cal T}^l={g^{li}\over 2}\left[2h^{\alpha\beta}{\partial g_{ij}\over\partial
t^\alpha}x^j_\beta+U^{(\alpha)}_{(i)j}x^j_\alpha+{\partial U^{(\alpha)}_{(i)}\over\partial t^\alpha}+U^
{(\alpha)}_{(i)}H^\gamma_{\alpha\gamma}-{\partial F\over\partial x^i}\right].
\end{equation}

The geometrical object ${\cal T}=({\cal T}^l)$ is a {\it
d-tensor field}  on $E=J^1(T,M)$. It follows that ${\cal T}$ can be written
as the $h$-trace of the d-tensor $\displaystyle{T^{(l)}_{(\alpha)\beta}=
{h_{\alpha\beta}\over p}{\cal T}^l}$, where $p=\dim T$, that is, ${\cal T}^l=
h^{\alpha\beta}T^{(l)}_{(\alpha)\beta}$. Of course, this writing is
not unique  but it is a natural extension of the case $\dim T=1$.

Finally, we conclude that the spatial $h$-spray ${\cal G}=({\cal G}^l)$  is
the $h$-trace of the spatial spray
\begin{equation}
G^{(l)}_{(\alpha)\beta}={1\over  2}\Gamma^l_{jk}x^j_\alpha x^k_\beta+
T^{(l)}_{(\alpha)\beta},
\end{equation}
that is, ${\cal G}^l=h^{\alpha\beta}G^{(l)}_{(\alpha)\beta}$.
\begin{th}
The extremals of the energy functional ${\cal E}_L$attached to a Kronecker
$h$-regular Lagrangian function $L$ on
$J^1(T,M)$ are harmonic maps \cite{12} of the spray $(H,G)$ defined by the
temporal components
$$
H^{(i)}_{(\alpha)\beta}=\left\{\begin{array}{ll}\medskip
\displaystyle{-{1\over 2}H^1_{11}(t)y^i,}&p=1\\\medskip
\displaystyle{-{1\over 2}H^\gamma_{\alpha\beta}x^i_\gamma,}&p\geq 2
\end{array}\right.
$$
and the local spatial components $G^{(i)}_{(\alpha)\beta}=$
$$
=\left\{\begin{array}{ll}
\displaystyle{{h_{11}g^{ik}\over 4}\left[{\partial^2L\over\partial x^j\partial
y^k}y^j-{\partial L\over\partial x^k}+{\partial^2L\over\partial t\partial y^k}+
{\partial L\over\partial x^k}H^1_{11}+2h^{11}H^1_{11}g_{kl}y^l\right]},&p=1
\\\medskip
\displaystyle{{1\over 2}\Gamma^i_{jk}x^j_\alpha x^k_\beta+T^{(i)}_{(\alpha)\beta},}
&p\geq 2,
\end{array}\right.
$$
where $p=\dim T$.
\end{th}
\addtocounter{defin}{1}
{\bf Definition \thedefin}  The spray $(H,G)$ constructed in the preceding theorem is
called the {\it canonical spray attached to the metrical multi-time Lagrange
space $ML^n_p$.}

In the sequel, using the canonical spray $(H,G)$ of the metrical multi-time
Lagrange space $ML^n_p$, one naturally induces \cite{12} a {\it nonlinear connection} $\Gamma$
on $E=J^1(T,M)$, defined by the temporal components 
\begin{equation}
M^{(i)}_{(\alpha)\beta}=2H^{(i)}_{(\alpha)\beta}=\left\{\begin{array}{ll}\medskip
-H^1_{11}y^i,&p=1\\\medskip
-H^\gamma_{\alpha\beta}x^i_\gamma,&p\geq 2,
\end{array}\right.
\end{equation}
and the spatial components
\begin{equation}
\mbox{\hspace*{10mm}}
N^{(i)}_{(\alpha)j}={\partial{\cal G}^i\over\partial x^j_\gamma}h_{\alpha
\gamma}=\left\{\begin{array}{ll}\medskip
\displaystyle{h_{11}{\partial{\cal G}^i\over\partial y^j}},&p=1\\\medskip
\displaystyle{\Gamma^i_{jk}x^k_\alpha+{g^{ik}\over 2}{\partial g_{jk}\over
\partial t^\alpha}+{g^{ik}\over 4}h_{\alpha\gamma}U^{(\gamma)}_{(k)j}},&p\geq 2,
\end{array}\right.
\end{equation}
where ${\cal G}^i=h^{\alpha\beta}G^{(i)}_{(\alpha)\beta}$. The nonlinear
connection $\Gamma=(M^{(i)}_{(\alpha)\beta},N^{(i)}_{(\alpha)j})$ is called
the {\it canonical nonlinear connection of the metrical multi-time Lagrange
space $ML^n_p$.}\medskip\\
\addtocounter{rem}{1}
{\bf Remarks \therem} i) Considering the particular case $(T,h)=(R,\delta)$,
the canonical nonlinear connection $\Gamma=(0,N^{(i)}_{(1)j})$  of  the
relativistic rheonomic Lagrange space $RL^n=(J^1(R,M),L)$ reduces to the canonical nonlinear
connection of the Lagrange space $L^n=(M,L)$.

ii) In the case of  an autonomous electrodynamics metrical multi-time Lagrange space (i. e.,
$g_{ij}(t^\gamma,x^k,x^k_\gamma)=g_{ij}(x^k)$), the generalized Christoffel
symbols $\Gamma^i_{jk}(t^\mu,x^m)$ of the metrical d-tensor $g_{ij}$ reduce
to the classical ones $\gamma^i_{jk}(x^m)$ and the canonical nonlinear
connection becomes $\Gamma=(M^{(i)}_{(\alpha)\beta},N^{(i)}_{(\alpha)j})$,
where
\begin{equation}
M^{(i)}_{(\alpha)\beta}=\left\{\begin{array}{ll}\medskip
-H^1_{11}y^i,&p=1\\\medskip
-H^\gamma_{\alpha\beta}x^i_\gamma,&p\geq 2,
\end{array}\right.
\end{equation}
\begin{equation}
N^{(i)}_{(\alpha)j}=\left\{\begin{array}{ll}\medskip
\displaystyle{\gamma^i_{jk}y^k+{g^{ik}\over 4}h_{11}U^{(1)}_{(k)j}},&p=1\\\medskip
\displaystyle{\gamma^i_{jk}x^k_\alpha+{g^{ik}\over 4}h_{\alpha\gamma}U^
{(\gamma)}_{(k)j}},&p\geq 2.
\end{array}\right.
\end{equation}

\section{Cartan canonical $h$-normal $\Gamma$-linear connection.
d-Torsions and d-curvatures}

\setcounter{equation}{0}
\hspace{5mm} Suppose that $J^1(T,M)$ is endowed with a nonlinear connection
$\Gamma$ defined
by the temporal components $M^{(i)}_{(\alpha)\beta}$ and the spatial
components $N^{(i)}_{(\alpha)j}$. Let $\displaystyle{\left\{{\delta\over\delta
t^\alpha}, {\delta\over\delta x^i}, {\partial\over\partial x^i_\alpha}\right\}
\subset{\cal X}(E)}$ and $\{dt^\alpha, dx^i, \delta x^i_\alpha\}\subset{\cal
X}^*(E)$ be the adapted bases of the nonlinear connection $\Gamma$, where
\begin{equation}
\left\{\begin{array}{l}\medskip
\displaystyle{{\delta\over\delta t^\alpha}={\partial\over\partial t^\alpha}-
M^{(j)}_{(\beta)\alpha}{\partial\over\partial x^j_\beta}}\\\medskip
\displaystyle{{\delta\over\delta x^i}={\partial\over\partial x^i}-
N^{(j)}_{(\beta)i}{\partial\over\partial x^j_\beta}}\\
\delta x^i_\alpha=dx^i_\alpha+M^{(i)}_{(\alpha)\beta}dt^\beta+N^{(i)}_{(\alpha)
j}dx^j.
\end{array}\right.
\end{equation}
Using the notations
$$
\begin{array}{lll}\medskip
\displaystyle{{\cal X}({\cal H}_T)=Span\left\{{\delta\over\delta t^\alpha}\right
\}},&
\displaystyle{{\cal X}({\cal H}_M)=Span\left\{{\delta\over\delta x^i}\right
\}},&
\displaystyle{{\cal X}({\cal V})=Span\left\{{\partial\over\partial x^i_\alpha}
\right\}},\\
{\cal X}^*({\cal H}_T)=Span\{dt^\alpha\},&
{\cal X}^*({\cal H}_M)=Span\{dx^i\},&
{\cal X}^*({\cal V})=Span\{\delta x^i_\alpha\},
\end{array}
$$
we obtain without difficulties the following
\begin{prop}
i) The Lie algebra  ${\cal X}(E)$ of vector fields decomposes as
$$
{\cal X}(E)={\cal X}({\cal H}_T)\oplus{\cal X}({\cal H}_M)\oplus{\cal X}(
{\cal V}).
$$

ii) The Lie algebra  ${\cal X}^*(E)$ of covector fields decomposes as
$$
{\cal X}^*(E)={\cal X}^*({\cal H}_T)\oplus{\cal X}^*({\cal H}_M)\oplus{\cal X}
^*({\cal V}).
$$
\end{prop}

Let us consider $h_T$, $h_M$ and $v$ the canonical projections of the above
decompositions.\medskip\\
\addtocounter{defin}{1}
{\bf Definition \thedefin} A linear connection $\nabla:{\cal X}(E)\times
{\cal X}(E)\to{\cal X}(E)$ is called a {\it $\Gamma$-linear connection} on
$E$ if $\nabla h_T=0,\;\nabla h_M=0,\;\nabla v=0$.\medskip

In the adapted basis
$\displaystyle{\left\{{\delta\over\delta
t^\alpha}, {\delta\over\delta x^i}, {\partial\over\partial x^i_\alpha}\right\}
\subset{\cal X}(E)}$, a $\Gamma$-linear connection $\nabla$ on $E$ is defined
by nine local coefficients,
\begin{equation}\label{lglc}
\hspace*{5mm}\nabla\Gamma=(\bar G^\alpha_{\beta\gamma},G^k_{i\gamma},G^{(k)(\beta)}_
{(\alpha)(i)\gamma},\bar L^\alpha_{\beta j},L^k_{ij},L^{(k)(\beta)}_{(\alpha)
(i)j},\bar C^{\alpha(\gamma)}_{\beta(j)},C^{k(\gamma)}_{i(j)},C^{(k)(\beta)
(\gamma)}_{(\alpha)(i)(j)}),
\end{equation}
introduced by
$$
\begin{array}{lll}\medskip
\displaystyle{\nabla_{\delta\over\delta t^\gamma}{\delta\over\delta t^\beta}=
\bar G^\alpha_{\beta\gamma}{\delta\over\delta t^\alpha},}&
\displaystyle{\nabla_{\delta\over\delta t^\gamma}{\delta\over\delta x^i}=
G^k_{i\gamma}{\delta\over\delta x^k},}&
\displaystyle{\nabla_{\delta\over\delta t^\gamma}{\partial\over\partial x^i
_\beta}=G^{(k)(\beta)}_{(\alpha)(i)\gamma}{\partial\over\partial x^k_\alpha}
,}\\\medskip
\displaystyle{\nabla_{\delta\over\delta x^j}{\delta\over\delta t^\beta}=
\bar L^\alpha_{\beta j}{\delta\over\delta t^\alpha},}&
\displaystyle{\nabla_{\delta\over\delta x^j}{\delta\over\delta x^i}=
L^k_{ij}{\delta\over\delta x^k},}&
\displaystyle{\nabla_{\delta\over\delta x^j}{\partial\over\partial x^i
_\beta}=L^{(k)(\beta)}_{(\alpha)(i)j}{\partial\over\partial x^k_\alpha}
,}\\
\displaystyle{\nabla_{\partial\over\partial x^j_\gamma}{\delta\over\delta
t^\beta}=\bar C^{\alpha(\gamma)}_{\beta(j)}{\delta\over\delta t^\alpha},}&
\displaystyle{\nabla_{\partial\over\partial x^j_\gamma}{\delta\over\delta
x^i}=C^{k(\gamma)}_{i(j)}{\delta\over\delta x^k},}&
\displaystyle{\nabla_{\partial\over\partial x^j_\gamma}{\partial\over\partial
x^i_\beta}=C^{(k)(\beta)(\gamma)}_{(\alpha)(i)(j)}{\partial\over\partial x^k_
\alpha}}.
\end{array}
$$
\addtocounter{rem}{1}
{\bf Remark \therem} The transformation rules of the above connection coefficients
are completely described in \cite{12}.\medskip\\
\addtocounter{ex}{1}
{\bf Example \theex} If $h_{\alpha\beta}$ (resp. $g_{ij}$) is a
semi-Riemannian metric on the temporal (resp. spatial) manifold $T$ (resp.
$M$), $H^\gamma_{\alpha\beta}$ (resp. $\gamma^k_{ij}$) are its Christoffel
symbols and $\Gamma_0=(M^{(i)}_{(\alpha)\beta}, N^{(i)}_{(\alpha)j})$,
where $M^{(i)}_{(\alpha)\beta}=-H^\gamma_{\alpha\beta}x^i_\gamma,\;N^{(i)}_
{(\alpha)j}=\gamma^i_{jk}x^k_\alpha$, is the canonical nonlinear connection
on $E$ attached to the metric pair $(h_{\alpha\beta},g_{ij})$, then
the following set of local coefficients \cite{12}
$$
B\Gamma_0=(\bar G^\alpha_{\beta\gamma},0,G^{(k)(\beta)}_{(\alpha)(i)\gamma},
0,L^k_{ij},L^{(k)(\beta)}_{(\alpha)(i)j},0,0,0),
$$
where $\bar G^\gamma_{\alpha\beta}=H^\gamma_{\alpha\beta},\;G^{(k)(\beta)}_
{(\gamma)(i)\alpha}=-\delta^k_iH^\beta_{\alpha\gamma},\;L^k_{ij}=\gamma_{ij}
^k$ and $L^{(k)(\beta)}_{(\gamma)(i)j}=\delta^\beta_\gamma\gamma^k_{ij}$, is
a\linebreak $\Gamma_0$-linear connection which is called {\it the Berwald
$\Gamma_0$-linear connection of the metric pair $(h_{\alpha\beta},g_{ij})$}.
\medskip

We recall that a $\Gamma$-linear connection $\nabla$ on $E$, defined by the local
coefficients \ref{lglc}
induces a natural linear connection on the d-tensors set of the jet fibre
bundle $E=J^1(T,M)$, which is characterized by a collection of local derivative
operators like "$ _{/\varepsilon}$", "$_{\vert p}$" and "$\vert
^{(\varepsilon)}_{(p)}$". The previous local operators are called the {\it
$T$-horizontal covariant derivative, $M$-horizontal covariant derivative}
and {\it vertical covariant derivative}. The detalied expressions of these
derivative operators are completely described in \cite{15}.

The study of the torsion {\bf T} and curvature {\bf R} d-tensors of an arbitrary
$\Gamma$-linear connection $\nabla$ on $E$ was made  in \cite{14}. In this
context, we proved that the torsion d-tensor is
determined by twelve effective local torsion d-tensors, while the curvature
d-tensor of $\nabla$ is determined by eighteen local d-tensors.

Now, let $h_{\alpha\beta}$ be a fixed pseudo-Riemannian metric on the temporal
manifold $T$, $H^\gamma_{\alpha\beta}$ its Christoffel symbols and $J=J^{(i)}
_{(\alpha)\beta j}{\partial\over\partial x^i_\alpha}\otimes dt^\beta\otimes
dx^j$, where $J^{(i)}_{(\alpha)\beta j}=h_{\alpha\beta}\delta^i_j$, the
normalization d-tensor \cite{12} attached to the metric $h_{\alpha\beta}$. The
big number of torsion and curvature d-tensors which characterize a general
$\Gamma$-linear connection on $E$ determines us to consider the following
\cite{11}\medskip\\
\addtocounter{defin}{1}
{\bf Definition  \thedefin} A $\Gamma$-linear connection $\nabla$ on $E=J^1(T,M)$,
defined by the local coefficients
$$
\nabla\Gamma=(\bar G^\alpha_{\beta\gamma},G^k_{i\gamma},G^{(k)(\beta)}_
{(\alpha)(i)\gamma},\bar L^\alpha_{\beta j},L^k_{ij},L^{(k)(\beta)}_{(\alpha)
(i)j},\bar C^{\alpha(\gamma)}_{\beta(j)},C^{k(\gamma)}_{i(j)},C^{(k)(\beta)
(\gamma)}_{(\alpha)(i)(j)}),
$$
that verify the relations $\bar G^\alpha_{\beta\gamma}=H^\alpha_{\beta\gamma},
\;\bar L^\alpha_{\beta j}=0,\;\bar C^{\alpha(\gamma)}_{\beta(j)}=0$ and
$\nabla J=0$, is called a {\it $h$-normal $\Gamma$-linear connection}.\medskip\\
\addtocounter{rem}{1}
{\bf Remark \therem} Taking into account the local covariant $T$-horizontal
"$_{/\gamma}$", $M$-horizontal "$_{\vert k}$" and vertical "$\vert^{(
\gamma)}_{(k)}$" covariant derivatives induced by $\nabla$, the condition
$\nabla J=0$ is equivalent to
$$
J^{(i)}_{(\alpha)\beta j/\gamma}=0,\quad
J^{(i)}_{(\alpha)\beta j\vert k}=0,\quad
J^{(i)}_{(\alpha)\beta j}\vert^{(\gamma)}_{(k)}=0.
$$

In this context, it is proved in \cite{11} the following
\begin{th}
The coefficients of a $h$-normal $\Gamma$-linear
connection $\nabla$ verify the identities
$$
\begin{array}{lll}\medskip
\bar G^\gamma_{\alpha\beta}=H^\gamma_{\alpha\beta},&\bar L^\alpha_{\beta j}=0,&
\bar C^{\alpha(\gamma)}_{\beta(j)}=0,\\\medskip
G^{(k)(\beta)}_{(\alpha)(i)\gamma}=\delta^\beta_\alpha G^k_{i\gamma}-\delta^k_i
H^\beta_{\alpha\gamma},&L^{(k)(\beta)}_{(\alpha)(i)j}=\delta^\beta_\alpha
L^k_{ij},&C^{(k)(\beta)(\gamma)}_{(\alpha)(i)(j)}=\delta^\beta_\alpha
C^{k(\gamma)}_{i(j)}.
\end{array}
$$
\end{th}
\addtocounter{rem}{1}
{\bf Remarks \therem} i) The preceding theorem implies that a $h$-normal
$\Gamma$-linear on $E$ is determined just by four effective coefficients
$$
\nabla\Gamma=(H^\gamma_{\alpha\beta},G^k_{i\gamma},L^k_{ij},C^{k(\gamma)}_
{i(j)}).
$$

ii) In the particular case $(T,H)=(R,\delta)$, a $\delta$-normal
$\Gamma$-linear connection identifies to the notion of $N$-linear connection
used in \cite{7}.\medskip\\
\addtocounter{ex}{1}
{\bf Example \theex} The canonical Berwald $\Gamma_0$-linear connection associated
to the metric pair $(h_{\alpha\beta},g_{ij})$ is a $h$-normal
$\Gamma_0$-linear connection defined by the local coefficients
$B\Gamma_0=(H^\gamma_{\alpha\beta},0,\gamma^k_{ij},0)$.\medskip

Note that, in the particular case of a $h$-normal $\Gamma$-linear connection
$\nabla$, the torsion d-tensors $\bar T^\mu_{\alpha\beta},\;\bar T^\mu_
{\alpha j}$ and $\bar P^{\mu(\beta)}_{\alpha(j)}$ vanish. Thus, the torsion
d-tensor {\bf T} of $\nabla$ is determined by the following nine local
d-tensors \cite{11}
\begin{equation}
\begin{tabular}{|c|c|c|c|}
\hline
&$h_T$&$h_M$&$v$\\
\hline
$h_Th_T$&0&0&$R^{(m)}_{(\mu)\alpha\beta}$\\
\hline
$h_Mh_T$&0&$T^m_{\alpha j}$&$R^{(m)}_{(\mu)\alpha j}$\\
\hline
$h_Mh_M$&0&$T^m_{ij}$&$R^{(m)}_{(\mu)ij}$\\
\hline
$vh_T$&0&0&$P^{(m)\;\;(\beta)}_{(\mu)\alpha(j)}$\\
\hline
$vh_M$&0&$P^{m(\beta)}_{i(j)}$&$P^{(m)\;(\beta)}_{(\mu)i(j)}$\\
\hline
$vv$&0&0&$S^{(m)(\alpha)(\beta)}_{(\mu)(i)(j)}$\\
\hline
\end{tabular}
\end{equation}
where
$$
\displaystyle{P^{(m)\;\;(\beta)}_{(\mu)\alpha(j)}={\partial M^{(m)}_{(\mu)
\alpha}\over\partial x^j_\beta}-\delta^\beta_\mu G^m_{j\alpha}+\delta^m_j
H^\beta_{\mu\alpha},}\quad
\displaystyle{P^{(m)\;\;(\beta)}_{(\mu)i(j)}={\partial N^{(m)}_{(\mu)i}\over
\partial x^j_\beta}-\delta^\beta_\mu L^m_{ji},}
$$
$$
\displaystyle{R^{(m)}_{(\mu)\alpha\beta}={\delta M^{(m)}_{(\mu)\alpha}\over
\delta t^\beta}-{\delta M^{(m)}_{(\mu)\beta}\over\delta t^\alpha},}\quad
\displaystyle{
R^{(m)}_{(\mu)\alpha j}={\delta M^{(m)}_{(\mu)\alpha}\over
\delta x^j}-{\delta N^{(m)}_{(\mu)j}\over\delta t^\alpha},}
$$
\medskip
$$
\displaystyle{
R^{(m)}_{(\mu)ij}={\delta N^{(m)}_{(\mu)i}\over\delta x^j}-{\delta N^{(m)}_
{(\mu)j}\over\delta x^i},}\quad
S^{(m)(\alpha)(\beta)}_{(\mu)(i)(j)}=\delta^\alpha_\mu C^{m(\beta)}_{i(j)}-
\delta^\beta_\mu C^{m(\alpha)}_{j(i)}\;,
$$
\medskip
$$
T^m_{\alpha j}=-G^m_{j\alpha}\;,\quad
T^m_{ij}=L^m_{ij}-L^m_{ji}\;,\quad
P^{m(\beta)}_{i(j)}=C^{m(\beta)}_{i(j)}.
$$
\addtocounter{rem}{1}
{\bf Remark \therem} For the Berwald $\Gamma_0$-linear connection associated to the
metrics $h_{\alpha\beta}$ and $g_{ij}$, all torsion d-tensors vanish,
except
$$
R^{(m)}_{(\mu)\alpha\beta}=-H^\gamma_{\mu\alpha\beta}x^m_\gamma,\quad
R^{(m)}_{(\mu)ij}=r^m_{ijl}x^l_\mu,
$$
where $H^\gamma_{\mu\alpha\beta}$ (resp. $r^m_{ijl}$) are the curvature tensors
of the metric $h_{\alpha\beta}$ (resp. $g_{ij}$).\medskip

The number of the effective curvature d-tensors of a $h$-normal $\Gamma$-linear
connection $\nabla$ reduces from eighteen to seven. The local d-tensors of
the curvature d-tensor {\bf R} of $\nabla$ are represented in
the table \cite{11}
\begin{equation}
\begin{tabular}{|c|c|c|c|}
\hline
&$h_T$&$h_M$&$v$\\
\hline
$h_Th_T$&$H^\alpha_{\eta\beta\gamma}$&$R^l_{i\beta\gamma}$&
$R^{(l)(\alpha)}_{(\eta)(i)\beta\gamma}=\delta^\alpha_\eta R^l_{i\beta\gamma}
+\delta^l_iH^\alpha_{\eta\beta\gamma}$\\
\hline
$h_Mh_T$&0&$R^l_{i\beta k}$&$R^{(l)(\alpha)}_{(\eta)(i)\beta k}=\delta^\alpha
_\eta R^l_{i\beta k}$\\
\hline
$h_Mh_M$&0&$R^l_{ijk}$&$R^{(l)(\alpha)}_{(\eta)(i)jk}=\delta^\alpha_\eta R^l_{
ijk}$\\
\hline
$vh_T$&0&$P^{l\;\;(\gamma)}_{i\beta(k)}$&$P^{(l)(\alpha)\;\;(\gamma)}_{(\eta)
(i)\beta(k)}=\delta^\alpha_\eta P^{l\;\;(\gamma)}_{i\beta(k)}$\\
\hline
$vh_M$&0&$P^{l\;(\gamma)}_{ij(k)}$&$P^{(l)(\alpha)\;(\gamma)}_{(\eta)
(i)j(k)}=\delta^\alpha_\eta P^{l\;(\gamma)}_{ij(k)}$\\
\hline
$vv$&0&$S^{l(\beta)(\gamma)}_{i(j)(k)}$&$S^{(l)(\alpha)(\beta)(\gamma)}_
{(\eta)(i)(j)(k)}=\delta^\alpha_\eta S^{l(\beta)(\gamma)}_{i(j)(k)}$\\
\hline
\end{tabular}
\end{equation}
where\medskip

$
\displaystyle{H^\alpha_{\eta\beta\gamma}={\partial H^\alpha_{\eta\beta}\over
\partial t^\gamma}-{\partial H^\alpha_{\eta\gamma}\over\partial t^\beta}+
H^\mu_{\eta\beta}H^\alpha_{\mu\gamma}-H^\mu_{\eta\gamma}H^\alpha_{\mu\beta},}
$\medskip

$
\displaystyle{R^l_{i\beta\gamma}={\delta G^l_{i\beta}\over\delta t^\gamma}-
{\delta G^l_{i\gamma}\over\delta t^\beta}+G^m_{i\beta}G^l_{m\gamma}-
G^m_{i\gamma}G^l_{m\beta}+C^{l(\mu)}_{i(m)}R^{(m)}_{(\mu)\beta\gamma},}
$\medskip

$
\displaystyle{R^l_{i\beta k}={\delta G^l_{i\beta}\over\delta x^k}-
{\delta L^l_{ik}\over\delta t^\beta}+G^m_{i\beta}L^l_{mk}-
L^m_{ik}G^l_{m\beta}+C^{l(\mu)}_{i(m)}R^{(m)}_{(\mu)\beta k},}
$\medskip

$
\displaystyle{R^l_{ijk}={\delta L^l_{ij}\over\delta x^k}-
{\delta L^l_{ik}\over\delta x^j}+L^m_{ij}L^l_{mk}-L^m_{ik}L^l_{mj}+
C^{l(\mu)}_{i(m)}R^{(m)}_{(\mu)jk},}
$\medskip

$
\displaystyle{P^{l\;\;(\gamma)}_{i\beta(k)}={\partial G^l_{i\beta}\over\partial
x^k_\gamma}-C^{l(\gamma)}_{i(k)/\beta}+C^{l(\mu)}_{i(m)}P^{(m)\;\;(\gamma)}_
{(\mu)\beta(k)},}
$\medskip

$
\displaystyle{P^{l\;(\gamma)}_{ij(k)}={\partial L^l_{ij}\over\partial
x^k_\gamma}-C^{l(\gamma)}_{i(k)\vert j}+C^{l(\mu)}_{i(m)}P^{(m)\;(\gamma)}_
{(\mu)j(k)},}
$\medskip

$
\displaystyle{S^{l(\beta)(\gamma)}_{i(j)(k)}={\partial C^{l(\beta)}_{i(j)}
\over\partial x^k_\gamma}-{\partial C^{l(\gamma)}_{i(k)}\over\partial x^j_
\beta}+C^{m(\beta)}_{i(j)}C^{l(\gamma)}_{m(k)}-
C^{m(\gamma)}_{i(k)}C^{l(\beta)}_{m(j)}.}
$\medskip\\
\addtocounter{rem}{1}
{\bf Remark \therem} In the case of the Berwald $\Gamma_0$-linear connection
associated to the metric pair $(h_{\alpha\beta},g_{ij})$, all curvature
d-tensors vanish, except $H^\delta_{\alpha\beta\gamma}$ and $R^l_{ijk}=r^l_
{ijk}$, where $r^l_{ijk}$ is the curvature tensor of the metric $g_{ij}$.
\medskip

Now, let  us consider $ML^n_p=(J^1(T,M),L)$ a metrical multi-time Lagrangian space and
$$
G^{(\alpha)(\beta)}_{(i)(j)}={1\over 2}{\partial^2L\over\partial x^i_\alpha
\partial x^j_\beta}=\left\{\begin{array}{ll}\medskip
h^{11}(t)g_{ij}(t,x^k,y^k),&p=1\\
h^{\alpha\beta}(t^\gamma)g_{ij}(t^\gamma,x^k),&p\geq 2
\end{array}\right.
$$
its vertical fundamental metrical d-tensor. Let $\Gamma=(M^{(i)}_{(\alpha)\beta},N^{(i)}
_{(\alpha)j})$ be the canonical nonlinear connection of the metrical
multi-time Lagrange space $ML^n_p$.

The main result of this paper  is the theorem of existence of the {\it
Cartan canonical $h$-normal linear connection} $C\Gamma$ which  allow the
natural subsequent development of the  {\it metrical multi-time Lagrange
theory of physical fields} \cite{10}.
\begin{th}{(existence and uniqueness of Cartan canonical connection)}\\
On the metrical multi-time Lagrange space $ML^n_p=(J^1(T,M),L)$ endowed with its canonical
nonlinear connection $\Gamma$ there is a unique $h$-normal $\Gamma$-linear
connection
$$
C\Gamma=(H^\gamma_{\alpha\beta},G^k_{j\gamma},L^i_{jk},C^{i(\gamma)}
_{j(k)})
$$
having the metrical properties\medskip

i) $g_{ij\vert k}=0,\quad g_{ij}\vert^{(\gamma)}_{(k)}=0$,\medskip

ii) $\displaystyle{ G^k_{j\gamma}={g^{ki}\over 2}{\delta g_{ij}\over\delta
t^\gamma},\quad L^k_{ij}=L^k_{ji},\quad
C^{i(\gamma)}_{j(k)}=C^{i(\gamma)}_{k(j)}}$.
\end{th}

{\bf Proof.} Let $C\Gamma=(\bar G^\gamma_{\alpha\beta},G^k_{j\gamma},L^i_{jk},
C^{i(\gamma)}_{j(k)})$ be h-normal $\Gamma$-linear connection whose
coefficients are defined by
$\displaystyle{\bar G^\gamma_{\alpha\beta}=H^\gamma_{\alpha\beta},\;
G^k_{j\gamma}={g^{ki}\over 2}{\delta g_{ij}\over\delta t^\gamma},}$ and
\begin{equation}
\begin{array}{l}\medskip
\displaystyle{L^i_{jk}={g^{im}\over 2}\left({\delta g_{jm}\over\delta x^k}+
{\delta g_{km}\over\delta x^j}-{\delta g_{jk}\over\delta x^m}\right),}\\
\displaystyle{C^{i(\gamma)}_{j(k)}={g^{im}\over 2}\left({\partial g_{jm}\over\partial
x^k_\gamma}+{\partial g_{km}\over\partial x^j_\gamma}-{\partial g_{jk}\over
\partial x^m_\gamma}\right)}.
\end{array}
\end{equation}
By computations, one easily verifies that  $C\Gamma$ satisfies the  conditions
{\it i} and {\it  ii}.

Conversely, let  us consider  a h-normal
$\Gamma$-linear connection
$$
\tilde C\Gamma=(\tilde{\bar G}^\gamma_{\alpha\beta},
\tilde G^k_{j\gamma},\tilde L^i_{jk},\tilde C^{i(\gamma)}_{j(k)})
$$
which satisfies {\it  i} and {\it ii}. It follows that
$$
\displaystyle{\tilde{\bar G}^\gamma_{\alpha\beta}=H^\gamma_
{\alpha\beta},\;\mbox{and}\;\tilde G^k_{j\gamma}={g^{ki}\over 2}{\delta g_{ij}\over
\delta t^\gamma}}.
$$

The condition $g_{ij\vert k}=0$ is equivalent to
$$
{\delta g_{ij}\over\delta x^k}=g_{mj}\tilde L^m_{ik}+g_{im}\tilde L^m_{jk}.
$$
Applying a Christoffel process to the indices $\{i,j,k\}$, we find
$$
\displaystyle{\tilde L^i_{jk}={g^{im}\over 2}\left({\delta g_{jm}\over\delta x^k}+
{\delta g_{km}\over\delta x^j}-{\delta g_{jk}\over\delta x^m}\right)}.
$$

By  analogy, using the relations $C^{i(\gamma)}_{j(k)}=C^{i(\gamma)}_{k(j)}$
and $g_{ij}\vert^{(\gamma)}_{(k)}=0$, following a Christoffel process applied
to the indices $\{i,j,k\}$, we obtain
$$
\displaystyle{\tilde C^{i(\gamma)}_{j(k)}={g^{im}\over 2}\left({\partial g_{jm}\over\partial
x^k_\gamma}+{\partial g_{km}\over\partial x^j_\gamma}-{\partial g_{jk}\over
\partial x^m_\gamma}\right)}.
$$

In conclusion, the uniqueness of the Cartan canonical  connection $C\Gamma$ is
clear. \rule{5pt}{5pt}\medskip\\
\addtocounter{rem}{1}
{\bf Remarks \therem} i) Replacing the canonical nonlinear connection $\Gamma$ by
a general one, the previous theorem holds good.

ii) In the particular case $(T,h)=(R,\delta)$, the Cartan canonical
$\delta$-normal $\Gamma$-linear connection of the relativistic rheonomic Lagrange space
$RL^n=(J^1(R,M),L)$ reduces to the Cartan canonical connection of the Lagrange
space $L^n=(M,L)$, constructed in \cite{7}.

iii) As a rule, the  Cartan canonical connection of a metrical multi-time Lagrange space
$ML^n_p$ verifies also the metrical properties
$$
h_{\alpha\beta/\gamma}=h_{\alpha\beta\vert k}=h_{\alpha\beta}\vert^{(\gamma)}
_{(k)}=0\;\mbox{and}\;g_{ij/\gamma}=0.
$$

iv) In the case
$\dim T\geq 2$, the coefficients of the Cartan connection of a metrical
multi-time Lagrange space reduce to
$$
\bar G^\gamma_{\alpha\beta}=H^\gamma_{\alpha\beta},\;G^k_{j\gamma}={g^{ki}
\over 2}{\partial g_{ij}\over\partial t^\gamma},\;L^i_{jk}=\Gamma^i_{jk},\;
C^{i(\gamma)}_{j(k)}=0.
$$

v) Particularly,  the coefficients of the Cartan connection of an autonomous
metrical multi-time Lagrange space of electrodynamics (i. e., $g_{ij}(t^\gamma,x^k,x^k_\gamma)
=g_{ij}(x^k)$) are the same with those of the Berwald connection, namely,
$C\Gamma=(H^\gamma_{\alpha\beta},0,\gamma^i_{jk},0)$. Note that the Cartan
connection is a $\Gamma$-linear connection, where $\Gamma$ is the canonical
nonlinear connection of the metrical multi-time Lagrangian space while the Berwald connection
is a $\Gamma_0$-linear connection, $\Gamma_0$ being the canonical nonlinear
connection  associated to the metric pair $(h_{\alpha\beta},g_{ij})$.
Consequently, the Cartan and Berwald connections are distinct.
\begin{th}
The torsion d-tensor {\bf T} of the Cartan canonical connection of a metrical
multi-time Lagrange space is determined by the local components
\begin{equation}
\begin{tabular}{|c|*{2}{c|}*{2}{c|}*{2}{c|}}
\hline
&\multicolumn{2}{|c}{$h_T$}&
\multicolumn{2}{|c}{$h_M$}&
\multicolumn{2}{|c|}{$v$}\\
\cline{2-7}
&$p=1$&$p\geq 2$&$p=1$&$p\geq 2$&$p=1$&$p\geq 2$\\
\hline
$h_Th_T$&0&0&0&0&0&$R^{(m)}_{(\mu)\alpha\beta}$\\
\hline
$h_Mh_T$&0&0&$T^m_{1j}$&$T^m_{\alpha j}$&$R^{(m)}_{(1)1j}$&
$R^{(m)}_{(\mu)\alpha j}$\\
\hline
$h_Mh_M$&0&0&0&0&$R^{(m)}_{(1)ij}$&$R^{(m)}_{(\mu)ij}$\\
\hline
$vh_T$&0&0&0&0&$P^{(m)\;\;(1)}_{(1)1(j)}$&$P^{(m)\;\;(\beta)}_{(\mu)\alpha(j)}$\\
\hline
$vh_M$&0&0&$P^{m(1)}_{i(j)}$&0&$P^{(m)\;(1)}_{(1)i(j)}$&0\\
\hline
$vv$&0&0&0&0&0&0\\
\hline
\end{tabular}
\end{equation}
where,

i) for $p=\dim T=1$, we have\medskip

\hspace*{25mm}$T^m_{1j}=-G^m_{j1}$, $\;P^{m(1)}_{i(j)}=C^{m(1)}_
{i(j)}$, $\;P^{(m)\;(1)}_{(1)1(j)}=-G^m_{j1}$,\medskip

\hspace*{25mm}$\displaystyle{P^{(m)\;(1)}_{(1)i(j)}={\partial N^{(m)}_{(1)i}\over\partial y^j}
-L^m_{ji}}$, $\quad\displaystyle{{\delta N^{(m)}_{(1)i}\over\delta x^j}-
{\delta N^{(m)}_{(1)j}\over\delta x^i}}$,\medskip

\hspace*{25mm}$\displaystyle{R^{(m)}_{(1)1j}=
-{\partial N^{(m)}_{(1)j}\over\partial t}+H^1_{11}\left[N^{(m)}_{(1)j}-{\partial
N^{(m)}_{(1)j}\over\partial y^k}y^k\right]}$;\medskip

ii) for $p=\dim T\geq 2$, denoting
$$
\begin{array}{l}\medskip
\displaystyle{F^m_{i(\mu)}={g^{mp}\over 2}\left[{\partial g_{pi}\over\partial
t^\mu}+{1\over 2}h_{\mu\beta}U^{(\beta)}_{(p)i}\right]},\\\medskip
\displaystyle{H^\gamma_{\mu\alpha\beta}={\partial H^\gamma_{\mu\alpha}\over
\partial t^\beta}-{\partial H^\gamma_{\mu\beta}\over\partial t^\alpha}+
H^\eta_{\mu\alpha}H^\gamma_{\eta\beta}-H^\eta_{\mu\beta}H^\gamma_{\eta\alpha},}\\
\displaystyle{r^m_{pij}={\partial \Gamma^m_{pi}\over\partial x^j}-
{\partial \Gamma^m_{pj}\over\partial x^i}+\Gamma^k_{pi}\Gamma^m_{kj}-\Gamma^
k_{pj}\Gamma^m_{ki},}
\end{array}
$$
we have
$$
\begin{array}{l}\medskip
T^m_{\alpha j}=-G^m_{j\alpha}$, $\;P^{m\;\;(\beta)}_{(\mu)\alpha
(j)}=-\delta^\beta_\gamma G^m_{j\alpha}$, $\;R^{(m)}_{(\mu)\alpha(j)}=-H^\gamma
_{\mu\alpha\beta}x^m_\gamma,\\\medskip
\displaystyle{R^{(m)}_{(\mu)\alpha j}=
-{\partial N^{(m)}_{(\mu)j}\over\partial t^\alpha}+{g^{mk}\over 2}H^\beta_
{\mu\alpha}\left[{\partial g_{jk}\over\partial t^\beta}+{h_{\beta\gamma}\over
2}U^{(\gamma)}_{(k)j}\right]},\\
\displaystyle{R^{(m)}_{(\mu)ij}=r^m_{ijk}x^k_\mu+\left[
F^m_{i(\mu)\vert j}-F^m_{j(\mu)\vert i}\right]}.
\end{array}
$$
\end{th}
\addtocounter{rem}{1}
{\bf Remark \therem} In the case of autonomous metrical multi-time Lagrange
space of electrodynamics
(i. e., $g_{ij}(t^\gamma,x^k,x^k_\gamma)=g_{ij}(x^k)$), all torsion d-tensors
of the Cartan connection vanish,  except
$$
R^{(m)}_{(\mu)\alpha\beta}=-H^\gamma_{\mu\alpha\beta}x^m_\gamma,\quad
R^{(m)}_{(\mu)\alpha j}=-{h_{\mu\eta}g^{mk}\over 4}\left[H^\eta_{\alpha\gamma}
U^{(\gamma)}_{(k)j}+{\partial U^{(\eta)}_{(k)j}\over\partial t^\alpha}\right],
$$
$$
R^{(m)}_{(\mu)ij}=r^m_{ijk}x^k_\mu+{h_{\mu\eta}g^{mk}\over 4}\left[
U^{(\eta)}_{(k)i\vert j}+U^{(\eta)}_{(k)j\vert i}\right],
$$
where $H^\gamma_{\mu\alpha\beta}$ (resp. $r^m_{ijk}$) are the curvature tensors
of the semi-Riemannian metric $h_{\alpha\beta}$ (resp. $g_{ij}$).
\begin{th}
The  curvature d-tensor {\bf R} of the Cartan
canonical connection is determined by the local components
$$
\begin{tabular}{|c|*{2}{c|}*{2}{c|}*{2}{c|}}
\hline
&\multicolumn{2}{|c}{$h_T$}&
\multicolumn{2}{|c}{$h_M$}&
\multicolumn{2}{|c|}{$v$}\\
\cline{2-7}
&$p=1$&$p\geq 2$&$p=1$&$p\geq 2$&$p=1$&$p\geq 2$\\
\hline
$h_Th_T$&0&$H^\alpha_{\eta\beta\gamma}$&0&$R^l_{i\beta\gamma}$&0&
$R^{(l)(\alpha)}_{(\eta)(i)\beta\gamma}$\\
\hline
$h_Mh_T$&0&0&$R^l_{i1k}$&$R^l_{i\beta k}$&$R^{(l)(1)}_{(1)(i)1k}=R^l_{i1k}$&
$R^{(l)(\alpha)}_{(\eta)(i)\beta k}$\\
\hline
$h_Mh_M$&0&0&$R^l_{ijk}$&$R^l_{ijk}$&$R^{(l)(1)}_{(1)(i)jk}=R^l_{ijk}$&
$R^{(l)(\alpha)}_{(\eta)(i)jk}$\\
\hline
$vh_T$&0&0&$P^{(l)\;\;(1)}_{i1(k)}$&0&$P^{(l)(1)\;(1)}_{(1)(i)1(k)}=
P^{(l)\;\;(1)}_{i1(k)}$&0\\
\hline
$vh_M$&0&0&$P^{l\;(1)}_{ij(k)}$&0&$P^{(l)(1)\;(1)}_{(1)(i)j(k)}=P^{l\;(1)}_
{ij(k)}$&0\\
\hline
$vv$&0&0&$S^{l(1)(1)}_{i(j)(k)}$&0&$S^{(l)(1)(1)(1)}_{(1)(i)(j)(k)}=S^{l(1)(1)}
_{i(j)(k)}$&0\\
\hline
\end{tabular}
$$
where
$R^{(l)(\alpha)}_{(\eta)(i)\beta\gamma}=\delta^\alpha_\eta R^l_{i\beta\gamma}+
\delta^l_i H^\alpha_{\eta\beta\gamma}$,
$R^{(l)(\alpha)}_{(\eta)(i)\beta k}=\delta^\alpha_\eta R^l_{i\beta k}$,
$R^{(l)(\alpha)}_{(\eta)(i)jk}=\delta^\alpha_\eta R^l_{ijk}$ and\medskip

i) for $p=\dim T=1$, we have\medskip

$
\displaystyle{R^l_{i1k}={\delta G^l_{i1}\over\delta x^k}-
{\delta L^l_{ik}\over\delta t}+G^m_{i1}L^l_{mk}-
L^m_{ik}G^l_{m1}+C^{l(1)}_{i(m)}R^{(m)}_{(1)1k},}
$\medskip

$
\displaystyle{R^l_{ijk}={\delta L^l_{ij}\over\delta x^k}-
{\delta L^l_{ik}\over\delta x^j}+L^m_{ij}L^l_{mk}-L^m_{ik}L^l_{mj}+
C^{l(1)}_{i(m)}R^{(m)}_{(1)jk},}
$\medskip

$
\displaystyle{P^{l\;(1)}_{i1(k)}={\partial G^l_{i1}\over\partial
y^k}-C^{l(1)}_{i(k)/1}+C^{l(1)}_{i(m)}P^{(m)\;(1)}_
{(1)1(k)},}
$\medskip

$
\displaystyle{P^{l\;(1)}_{ij(k)}={\partial L^l_{ij}\over\partial
y^k}-C^{l(1)}_{i(k)\vert j}+C^{l(1)}_{i(m)}P^{(m)\;(1)}_
{(1)j(k)},}
$\medskip

$
\displaystyle{S^{l(1)(1)}_{i(j)(k)}={\partial C^{l(1)}_{i(j)}
\over\partial y^k}-{\partial C^{l(1)}_{i(k)}\over\partial y^j}
+C^{m(1)}_{i(j)}C^{l(1)}_{m(k)}-C^{m(1)}_{i(k)}C^{l(1)}_{m(j)};}
$\medskip

ii) for $p=\dim T\geq 2$, we have\medskip

$
\displaystyle{H^\alpha_{\eta\beta\gamma}={\partial H^\alpha_{\eta\beta}\over
\partial t^\gamma}-{\partial H^\alpha_{\eta\gamma}\over\partial t^\beta}+
H^\mu_{\eta\beta}H^\alpha_{\mu\gamma}-H^\mu_{\eta\gamma}H^\alpha_{\mu\beta},}
$\medskip

$
\displaystyle{R^l_{i\beta\gamma}={\delta G^l_{i\beta}\over\delta t^\gamma}-
{\delta G^l_{i\gamma}\over\delta t^\beta}+G^m_{i\beta}G^l_{m\gamma}-G^m_
{i\gamma}G^l_{m\beta},}
$\medskip

$
\displaystyle{R^l_{i\beta k}={\delta G^l_{i\beta}\over\delta x^k}-
{\delta\Gamma^l_{ik}\over\delta t^\beta}+G^m_{i\beta}\Gamma^l_{mk}-\Gamma^m_
{ik}G^l_{m\beta},}
$\medskip

$
\displaystyle{R^l_{ijk}=r^l_{ijk}={\partial\Gamma^l_{ij}\over\partial x^k}-
{\partial\Gamma^l_{ik}\over\partial x^j}+\Gamma^m_{ij}\Gamma^l_{mk}-\Gamma^
m_{ik}\Gamma^l_{mj}}
$.
\end{th}
\addtocounter{rem}{1}
{\bf Remark \therem} In the case of an autonomous metrical multi-time Lagrange
space of electrodynamics
(i. e. , $g_{ij}(t^\gamma,x^k,x^k_\gamma)=g_{ij}(x^k)$), all curvature d-tensors
of the Cartan canonical connection vanish,  except $H^\alpha_{\eta\beta\gamma}$
and\linebreak
$R^l_{ijk}=r^l_{ijk}$, that is, the curvature tensors  of the semi-Riemannian
metrics $h_{\alpha\beta}$ and  $g_{ij}$.\medskip\\
{\bf\underline{Open problem}.} The development of an analogous metrical
multi-time Lagrangian geometry on $J^2(T,M)$ is in our attention.\medskip\\
{\bf Acknowledgements.} Many thanks  go to Prof. Dr. M. Matsumoto and Prof.
Dr. D. Opri\c s, the first  readers who advice us to distroy the initial version.

\begin{center}
University POLITEHNICA of Bucharest\\
Department of Mathematics I\\
Splaiul Independentei 313\\
77206 Bucharest, Romania\\
e-mail: mircea@mathem.pub.ro\\
e-mail: udriste@mathem.pub.ro
\end{center}

\end{document}